# The converse to Curtiss' theorem for one-sided moment generating functions

Patrick Chareka

Department of Mathematics, Statistics and Computer Science, St. Francis Xavier University, P.O. Box 5000, Antigonish, Nova Scotia, B2G 2W5, Canada. Email: pchareka@ns.sympatico.ca.

**Abstract**

A set of necessary and sufficient conditions for a sequence of moment generating functions to converge to a moment generating function on an interval (a,b) not necessarily containing 0, is given. The result is derived using recent results by Mukherjea, et al. (2006) and Chareka (2007).

## 1 Introduction

The moment generating function (mgf), of a random variable $X$ with distribution function $F(x) = P(X \leq x)$, is defined as

$$M(t) = \mathbb{E}\left(e^{tX}\right) = \int_{-\infty}^{\infty} e^{tx} dF(x). \tag{1}$$

where the integral on the right of equation (1) is the Riemann-Stieltjes integral of $e^{tx}$ with respect to $F(x)$. The domain of $M(t)$ is the set of all real $t$ for which the expectation $\mathbb{E}\left(e^{tX}\right)$ exists finitely. The mgf exists for $t = 0$. It is however, customary to say that the mgf exists if there exists a positive number $\delta$ such that $M(t)$ exists for all $t \in (-\delta, \delta)$. In this case, $X$ has has finite moments of all orders. It is possible for a random variable to have finite moments of all orders but when the corresponding mgf exists only for $t < 0$ (or $t > 0$). For example, the lognormal random variable $X = e^Z$, where $Z$ is a standard normal random variable has moment generating function

$$M(t) = \frac{1}{\sqrt{2\pi}} \int_0^\infty \frac{1}{x} \exp\left(tx - \frac{1}{2}(\ln x)^2\right) dx. \qquad (2)$$

The mgf in (2) exists for $t \leq 0$ but is infinite for $t > 0$. All moments of the lognormal distribution exist (finitely). It is also possible, for an mgf to exist for all $t$ in an interval not containing zero but when all the moments of the the corresponding distribution are infinite. For example, the Frechet distribution

$$F(x) = \begin{cases} 0, & x \leq 0 \\ \exp\left(-\frac{1}{x}\right), & x > 0, \end{cases} \qquad (3)$$

has moment generating function, $M(t) = \int_0^\infty \frac{1}{x^2} \exp\left(-\frac{1}{x} + tx\right) dx$. The moments of the Frechet distribution are all infinite, i.e. they do not exist. However, the corresponding mgf is finite for all $t \in (-\infty, 0]$ and infinite for $t > 0$. we use this example later in the paper.

The moment generating function has many theoretical and practical applications in probability and statistics. Most of the applications require that the mgf exist in some open interval containing zero. Examples of such results include the uniqueness theorem for moment generating functions, properties of locally sub-Gaussian random variables [2], and Curtiss' theorem for sequences of moment generating functions [4]. The famous Curtiss' theorem states that if $\{M_n(t)\}$ is a sequence of mgfs corresponding to a sequence of distribution functions $\{F_n(x)\}$, then convergence of $\{M_n(t)\}$ to a moment generating function $M(t)$ on $(-\delta, \delta)$ implies that $\{F_n(x)\}$ converges weakly to $F(x)$, where $F(x)$ is the distribution function with mgf $M(t)$.

It was shown recently in [7] that Curtiss' theorem does not require the open interval to include zero. More specifically it was shown that if a sequence of mgfs $\{M_n(t)\}$ converges pointwise to a moment generating function $M(t)$ for all $t$ in an open interval not necessarily containing the origin, then the corresponding distribution functions $\{F_n(x)\}$ converges weakly to the distribution $F(x)$ corresponding to $M(t)$. The essence of the result, is that the result tends to simplify proofs of certain limit theorems such as the central limit theorem based on Curtiss' theorem. It is clear that the result is weaker than Curtiss' theorem and generally easier to apply than Curtiss' theorem.

Following Curtiss' article on moment generating functions was another article by Kozakiewicz in 1947, [6]. The main motivation of Kozakiewicz's paper and this article, is the example in [4] which shows that in general, weak convergence of a sequence of distribution functions does not say much about the behavior of the corresponding sequence of moment generating functions. Kozakiewicz presents necessary and sufficient conditions for a sequence of moment generating functions to converge to a moment generating function in an open interval containing zero. Recently a new version of Curtiss' theorem for one-sided moment generating functions, i.e. moment generating functions that exist in an open interval not necessarily containing zero, was presented in [7]. The new theorem in [7], does not give the converse to Curtiss' theorem. In this article we present necessary and sufficient conditions for a sequence of moment generating functions to converge to a moment generating function in an open interval (not necessarily containing zero).

## 2 The converse to Curtiss' theorem for one-sided mgfs

**Theorem 1** . *Let $F_n(x)$ be the distribution function of a random variable $X_n$ and $M_n(t) = \int_{-\infty}^{\infty} e^{tx} dF_n(x)$ be the corresponding moment generating function, where $n \in \mathbb{N}$. Suppose that for each n, $M_n(t)$ exists for all $t \in (a,b)$. Then a set of necessary and sufficient conditions for the sequence $M_n(t)$ to converge to a moment generating function $M(t)$ on $(a,b)$, is,*

  (a) $\sup_{n \to \infty} M_n(t) < \infty, \ t \in (a,b)$,

  (b) $F_n(x)$ *converges weakly to $F(x)$ such that $\int_{-\infty}^{\infty} e^{tx} dF(x) = M(t)$ exists for all $t \in (a,b)$.*

**Proof:** It suffices to prove the result for $t \in (0,b)$. The proof for the case $t < 0$ is similar. First we note that existence of an mgf on any interval (a,b) uniquely determines a probability distribution (see Chareka, 2007).

*Necessity*: Suppose that $\lim_{n \to \infty} M_n(t) = M(t)$ for all $t \in (0,b)$, where $M(t)$ is the moment generating function of a random variable $X$. Then, it is clear that for each $t \in (0,b)$, the sequence $\{M_n\}$ is bounded, since every convergent sequence is bounded. Hence condition (a) of theorem (1) holds. Condition (b) follows from theorem 2 given in [7].

*Sufficiency*: Conversely suppose that conditions (a) and (b) hold. Then for each $t \in (0,b)$, $e^{tX_n}$ converges weakly to $e^{tX}$, since $e^{tX_n}$ is a continuous transformation of $X_n$. It follows from condition (a) and theorem 4.5.2 in [3] that $\lim_{n \to \infty} \mathbb{E}\left(e^{tX_n}\right) = \mathbb{E}\left(e^{tX}\right) = M(t)$ for all $t \in (0,b)$. This completes the proof of theorem (1).

An example to illustrate Curtiss' theorem for one-sided mgfs is given in [7]. Here we give an example to illustrate the converse of Curtiss' theorem for one-sided mgfs. Let $X_n$ have the distribution function,

$$F_n(x) = \begin{cases} 0, & x \leq 1/n \\ \left(1 - \frac{1}{nx}\right)^n, & x > 1/n. \end{cases} \quad (4)$$

When $n = 1$, $F_1(x)$ is the distribution function of a Pareto distribution with support $[1, \infty)$. Using the fact that $\left(1 - \frac{1}{nx}\right)^n \to e^{-1/x}$ for all $x \neq 0$, it is easy to see that $F_n(x)$ converges to the standard Frechet distribution given in (3). That is, condition (b) of theorem 1 is satisfied. Furthermore, it may be readily verified that for any $a < 0$, the mgf corresponding to $F_n(x)$ is bounded for all $t \in (a, 0]$. That is condition (b) of theorem (1) is satisfied. It is also easy to see that as $n \to \infty$,

$$M_n(t) = \int_0^\infty \frac{1}{x^2}\left(1 - \frac{1}{nx}\right)^{n-1} e^{tx} dx \to \int_0^\infty \frac{1}{x^2} \exp\left(-\frac{1}{x} + tx\right) dx. \quad (5)$$

The result in equation (5) may be deduced quickly from theorem (1) or proved by applying the dominated convegence theorem.

Condition (a) of Kozakiewicz's third theorem, is in general, not easy to verify. From the uniqueness theorem for mgfs given in Chareka (2007), condition (a) of theorem (1) in this paper is relatively much easier to verify (for some convenient short interval $(a, b)$), than Kozakiewicz's condition for the converse to Curtiss' theorem. An important and interesting case in which it may not even be necessary to verify condition (a) of theorem (1) is given in theorem (2) below.

**Theorem 2**. *Let $\{F_n(x)\}$ be a sequence of continuous distribution functions such that $M_n(t) = \int_0^\infty e^{tx} dF_n(x)$ exists for all $t \in (a, b)$. Suppose also that $F_n(x)$ converges weakly to a continuous distribution function $F(x)$ with mgf $M(t)$ which exists for all $t \in (a, b)$. Then $\lim_{n\to\infty} M_n(t) = M(t)$ for all $t \in (a, b)$.*

**Proof:** It is well-known, see for example, [8] that if $h_n(x)$ is a sequence of continuous and integrable functions converging uniformly to an integrable function $h(x)$ on a closed interval $I$, not necessarily finite, for example $I = [0, \infty)$, then $h(x)$ is continuous and

$$\lim_{n \to \infty} \int_0^\infty h_n(x) dx = \int_0^\infty h(x) dx. \quad (6)$$

It is important to note here that equation (6) implies the existence of the limit on the left-hand side. Now let $G_n((t, x) = P(e^{tX_n} \leq x) = P(X_n \leq \ln(x)/t) =$ $F_n(\ln(x)/t)$, $t \neq 0$, and $h_n(x) = 1 - G_n(t,x)$. Since for each $t \in (a,b)$, $e^{tX_n}$ is a continuous transformation of $X_n$, it follows that $G_n(t,x)$ converges weakly to $G(t,x) = F(\ln(x)/t)$. It is also known, see for example [5], that if a sequence of distribution functions converges to a continuous distribution function then the convergence must be uniform. Since $G(t,x)$ is continuous in $x$, it follows that for each $t$, $G_n(t,x)$ converges uniformly to $G(t,x)$. Hence,

$$\begin{align}
\lim_{n \to \infty} M_n(t) &= \lim_{n \to \infty} \mathbb{E}\left(e^{tX_n}\right) \tag{7}\\
&= \lim_{n \to \infty} \int_0^\infty (1 - G_n(t,x))\, dx \tag{8}\\
&= \int_0^\infty (1 - G(t,x))\, dx \tag{9}\\
&= \mathbb{E}\left(e^{tX}\right) \tag{10}\\
&= M(t). \tag{11}
\end{align}$$

For a sequence of independent and identically distributed continuous random variables $\{X_n\}$, each having mean $\mu$ and variance $\sigma^2$, an example of a new sequence of random variables converging to a continuous distribution is $\left\{Y_n = \frac{\overline{X} - \mu}{\sigma/\sqrt{n}}\right\}$, where $\overline{X} = (1/n)\sum_{i=1}^n X_i$. By the central limit theorem, $Y_n$ converges in distribution, to a standard normal random variable. It follows from theorem (2) that if $M_n(t) = \mathbb{E}\left(e^{tY_n}\right)$ exists for all $t \in (a,b)$, then $\lim_{n \to \infty} M_n(t) = \exp\left(-t^2/2\right)$. This is another example of the converse to Curtiss' theorem.

## 3 Conclusion

In this article we have derived the converse to Curtiss' theorem for one-sided moment generating functions. The result completes the work of Mukherjea, Rao and Suen and generalizes Kozakiewicz's necessary and sufficient conditions for the convergence of moment generating functions in an open interval containing zero, to convergence in an open interval not necessarily containing zero.